\documentclass[11pt, a4paper]{amsart}
\usepackage{fullpage} 
\usepackage{amsmath}
\usepackage{amsthm}
\usepackage{amssymb}
\usepackage{xypic}
\usepackage[only,mapsfrom]{stmaryrd}
\usepackage{graphicx}
\usepackage{cancel} 
\usepackage{color} 
\usepackage[all]{xy} 
\newtheorem{theorem}{Theorem}[section]
\newtheorem{lemma}[theorem]{Lemma}

\newtheorem{prop}[theorem]{Proposition}

\theoremstyle{definition}
\newtheorem{definition}[theorem]{Definition}
\newtheorem{example}[theorem]{Example}

\theoremstyle{remark}

\numberwithin{equation}{section}

\newcommand{\C}{\mathbb{C}}
\newcommand{\N}{\mathbb{N}}

\newcommand{\Z}{\mathbb{Z}}

\DeclareMathOperator{\Hom}{Hom}

\DeclareMathOperator{\Id}{Id}
\DeclareMathOperator{\id}{id} 
\DeclareMathOperator{\MU}{MU}
\DeclareMathOperator{\Ann}{Ann}
\title{An equivariant Quillen theorem}
\dedicatory{Dedicated to Tammo tom Dieck on the occasion of his eightieth birthday}

\author{Bernhard Hanke}   
\address{Institut f\"ur Mathematik, Universit\"at Augsburg, D-86135 Augsburg, Germany} 
\email{ hanke@math.uni-augsburg.de}

\author{Michael Wiemeler} 

\address{Mathematisches Institut, WWU M\"unster, Einsteinstr. 62, D-48149 M\"unster, Germany} 
\email{ wiemelerm@uni-muenster.de}



\date{\today }

\begin{document}

\begin{abstract}
A classical theorem due to Quillen (1969) identifies the unitary bordism ring with the Lazard ring, which classifies
the universal one-dimensional commutative formal group law. We prove an equivariant generalization 
of this result by identifying the homotopy theoretic $\Z/2$-equivariant unitary bordism ring, introduced 
by tom Dieck (1970), with the $\Z/2$-equivariant Lazard ring, introduced by Cole-Greenlees-Kriz (2000). 
Our proof combines a computation of the homotopy theoretic $\Z/2$-equivariant unitary bordism ring due to 
Strickland (2001) with a detailed investigation of the $\Z/2$-equivariant Lazard ring. \end{abstract}

\maketitle

\section{Introduction}

Almost 50 years ago Daniel Quillen gave the following algebraic description of the unitary bordism ring.

\begin{theorem}[\cite{Quillen}]  The  canonical map 
\[
     L _*  \to \MU_{*}  
\]
from the Lazard ring  to the coefficient ring of unitary bordism theory is an isomorphism. 
\end{theorem} 

Recall that the Lazard ring is the representing ring of the universal one-dimensional commutative formal group law 
\cite{Laz55}, and that the unitary bordism ring $\MU_* = \MU^{-*}$ is the underlying ring of a  one-dimensional formal group law
\[
     F(x,y) \in \MU^{-*}[[x,y]] = \MU^{-*}(\C P^{\infty} \times \C P^{\infty}) 
\]
given by the pull back of the universal Chern class in $\MU^2(\C P^{\infty})$ under the classifying map 
$\C P^{\infty} \times \C P^{\infty} \to \C P^{\infty}$ of the tensor product of the universal line bundles on each 
factor of $\C P^{\infty} \times \C P^{\infty}$.

Today Quillen's theorem  is  one of the organizational principles of stable homotopy theory. 
Establishing equivariant analogues of this result  is therefore a reasonable goal. 
Tom Dieck \cite{tD70} defined homotopy theoretic $G$-equivariant unitary bordism 
theories $\MU^{G}_*$ for compact Lie groups $G$, using 
arbitrary unitary $G$-representation as suspension coordinates.  Further information on the 
foundations of equivariant stable homotopy theory  is contained in  
\cite{LMS85}, for instance. 

For compact abelian Lie groups $A$ the notion of $A$-equivariant formal 
group laws was introduced in \cite{CGK00}, and subject to an extensive theoretical investigation in \cite{Str11}. 
Similarly to the non-equivariant situation there is a universal one-dimensional commutative $A$-equivariant formal 
group law, see \cite[Cor. 14.3]{CGK00}, together with a representing ring 
$L_A$,  and a classifying map
\[
  \lambda_A :   L_A \to \MU^A_* . 
 \]
An equivariant version of Quillen's theorem amounts to $\lambda_A$ being an isomorphism. 

Let $A$ be a finite abelian group. By use of a  localization-completion pull back square due to tom Dieck (for cyclic $A$ this is 
 \cite[Theorem 5.1]{tD70}), 
together with the classification of Euler-complete and Euler local equivariant 
formal group laws, one can show that $\lambda_A$ is surjective 
with each element in the kernel  being  Euler torsion and infinitely Euler divisible, see  \cite[Theorem 13.1]{G01}.

Strickland \cite{Str01} presented an algebraic 
description of $\MU^{\Z/2}_*$ in terms of generators and relations, 
and stated without proof the existence of a section $\MU^{\Z/2}_* \to L_{\Z/2}$ of the classifying map $\lambda_{\Z/2}$, 
establishing  $\MU_*^{\Z/2}$ as a retract of $L_{\Z/2}$. 
 
In spite of these positive results, injectivity of $\lambda_A$ has remained  elusive for any non-trivial $A$. 
This problem has been raised at several places, see e.g.~\cite[Questions 16.8]{G01} and  \cite[Chapter 13]{Str11}. 
We remark that there do exist non-additive $\Z/2$-equivariant formal group laws with representing rings 
containing non-zero, infinitely Euler divisible and Euler torsion
elements, see Example \ref{exotic}.  Hence injectivity of $\lambda_A$ does not 
follow merely from the results mentioned before, but 
requires some new structural insight concerning the ring $L_A$ itself.  This is what 
we will achieve in the paper at hand for the simplest non-trivial case $A = \Z/2$, providing 
 the first instance 
of an equivariant Quillen theorem:

\begin{theorem} \label{main}
The map  $\lambda_{\Z/2} :  L_{\Z/2} \rightarrow \MU^{\Z/2}_*$ is an isomorphism. 
\end{theorem}

Our argument starts with the construction of an explicit section $\mu_{\Z/2}$ of the classifying map 
$\lambda_{\Z/2}$. This 
allows us to  introduce structure constants $\rho_{ij} \in L_{\Z/2}$ in 
the  kernel of $\lambda_{\Z/2}$, which measure the deviation from $\lambda_{\Z/2}$ being 
an isomorphism. 

The proof of the vanishing of 
all $\rho_{ij}$ rests  on two major lines of 
thought, developed in Sections \ref{Eulerone} and \ref{new} of our paper. 
 The first one is of a conceptual nature and based on 
the construction of a {\em normalization functor},  turning $\Z/2$-equivariant formal 
group laws into ones with Euler classes equal to $1$. Applied 
to the universal $\Z/2$-equivariant formal group law we can hence 
derive efficient upper bounds on  the Euler torsion of the elements $\rho_{ij}$, see Theorem \ref{step1}. 

The second one is an explicit computation of a particular $\Z/2$-equivariant 
formal group law, which is obtained from $L_{\Z/2}$ by dividing out the 
ideal $\overline{J}$ generated by the images under $\mu_{\Z/2}$ of the positive degree generators of $\MU^{\Z/2}_*$. 
It turns out that the resulting $\Z/2$-equivariant 
formal group law is the additive one.  In other words (Theorem \ref{step2}): The structure constants $\rho_{ij}$ 
lie in $\overline{J}$.

Combining these two results in a bootstrap like manner 
forces vanishing of the $\rho_{ij}$. See Section \ref{strategy} 
for more details. 

On the one hand  we are optimistic that our approach can 
be generalized to more general $A$, for example 
cyclic groups of prime order. On the other hand we feel that 
the proof of an equivariant Quillen theorem for all compact abelian $A$ 
requires  some additional insight, which, among others, avoids explicit computations of $\MU^{A}_*$ 
in terms of generators and relations. We leave this topic for future research. 

\bigskip 

{\em Acknowledgements:}  We thank  the referee for a number of helpful comments. B.H. is grateful to the MPI,  
Bonn, to the IMPA, Rio de Janeiro, and to the IHES, Bures-sur-Yvette, 
for their hospitality while parts of this research were carried out.
This work was partially supported by a DFG research grant (B.H. and M.W.), by the DFG Priority 
Programme SPP 2026 (B.H. and M.W.), and by the SFB 878 at WWU M\"unster (M.W.). 

\section{Recollections on equivariant formal group laws} \label{recoll} 

For the notion of equivariant formal group laws, and the basics of the corresponding theory we refer the reader to \cite{CGK00}. 
Here we only recall some of the most important features and fix some notation. 
Let $A$ be a compact abelian Lie group. We denote by $A^* = \Hom(A, S^1)$ the group of irreducible unitary  $A$-representations. 
The trivial one-dimensional representation  is denoted 
by $\epsilon$. An {\em $A$-equivariant formal group law} is given by a 
quintuple  $(k, R, \Delta, \theta, y(\epsilon))$ with a commutative ring $k$, a complete commutative 
topological $k$-algebra $R$, a continuous comultiplication 
\[
   \Delta: R \to R \widehat{\otimes} R ,
\]
an augmentation $\theta: R \to k^{A^*}$ and an orientation $y(\epsilon) \in R$, such that the following
axioms are satisfied, see \cite[Def. 11.1]{CGK00}, which we here only recall for finite abelian $A$.

\begin{enumerate}
\item The comultiplication \(\Delta\)  is a map of \(k\)-algebras, co-commutative, co-associative and co-unital.
\item The augmentation \(\theta\) is a map of \(k\)-algebras compatible with the coproduct, so that \(\ker \theta\) defines the topology.
\item \(y(\epsilon)\in R\) is a regular element in the kernel of $\theta(\epsilon)$, and \(\theta(\epsilon)\) induces an isomorphism \(R/(y(\epsilon))\cong k\).
\end{enumerate}

We obtain an action $l$ of $A^*$ on $R$ by the formula
\[
    l_{\alpha} (c) = (\theta(\alpha^{-1}) \otimes \id) (\Delta(c))
\]
for $\alpha \in A^*$ and $c \in R$. Moreover,  corresponding to \(\alpha\in A^*\) we have coordinates  $y(\alpha) := l_{\alpha}(y(\epsilon)) \in R$. Finally, we define Euler 
classes 
\[
   e(\alpha) := \theta(\epsilon) ( y(\alpha))  \in k. 
\]
For explicit computations it is necessary to choose a basis of the topological 
$k$-module $R$. Recall \cite[Notation 12.1]{CGK00} 
that a  {\em complete flag} $F = (V^0 \subset V^1 \subset V^2 \subset 
\ldots)$  is a sequence of  \(r\)-dimensional  complex \(A\)-representations $V^r$ such that \(V^r\subset V^{r+1}\) and each finite dimensional complex \(A\)-representation is isomorphic to a subrepresentation of some \(V^r\). 
Given a complete flag $(V^r)_{r \in \N}$ we 
obtain a \(k\)-basis $(y(V^r))_{r \in \N}$ of the topological $k$-module $R$, see \cite[Lemma  13.2]{CGK00}. 
In this basis the coproduct $\Delta$ is given by 
\[
   \Delta(y(V^r)) = \sum_{i,j \geq 0} \beta^{(r)}_{i,j} \cdot y(V^i) \otimes y(V^j) 
\]
with structure constants $\beta^{(r)}_{i,j} \in k$. Let $k' \subset k$ be the subring 
generated by  the coefficients $\beta^{(1)}_{i,j}$ and the Euler classes $e(\alpha)$, $\alpha\in A^*$, and 
let $R' \subset R$ be the free $k'$-module with basis $(y(V^r))_{r\in \N}$. 
By an argument similar as for the proof of \cite[Theorem 16.1]{CGK00} the coproduct 
$\Delta$ restricts to a coproduct on $R'$ and thus we obtain 
an induced $\Z/2$-equivariant formal group law $(k', R', \Delta, \theta, y(\epsilon))$. 
We can hence assume without loss of generality that the underlying ring $k$ of a formal 
group law $(k,R)$ is generated by the structure constants 
$\beta^{(1)}_{i,j} \in k$ and the Euler classes $e(\alpha)$. In this case we briefly call $k$ the 
{\em representing ring} of the $A$-equivariant formal group law.

It is sometimes important to  consider {\em graded} $A$-equivariant formal group laws  $(k, R)$. 
This means that $k$ and $R$ 
are $\Z$-graded rings, where for  $r \geq 0$  the basis element \(y(V^r)\) sits in degree \(-2r\), and the coproduct, the 
$A^*$-action and the augmentation are grading preserving. In this case 
the  structure constants  \(\beta^{(r)}_{i,j}\) are homogeneous of degree \(2(i+j-r)\) and   the Euler classes 
have degree \(-2\). The most prominent example  is the universal $A$-equivariant formal group law 
$(L_A, R)$. Here we notice that the  construction of $(L_A, R)$ in  \cite[Cor. 14.3]{CGK00} in fact produces 
a graded $A$-equivariant formal group law. With this grading the classifying map 
\[
    L_A \to \MU^A_* 
\]
of the formal group law associated with $A$-equivariant unitary  bordism theory is 
grading preserving.  We remark that  this grading  structure of $L_A$, 
which plays an important role for our argument,  is not 
present in \cite{CGK00}. 

\section{Coordinate change} \label{sec:base_change} 

Let $A = \Z/2$ and let $(k,R,\Delta, \theta, y(\epsilon))$ be an $A$-equivariant formal group law.
We will work out some explicit formulas relating expansions of elements in $R$ with respect to different complete 
flags. Let $\eta \in A^*$ be the unique 
non-trivial one dimensional unitary $A$-representation and define
$e := e(\eta) \in k$ as the corresponding Euler class. For $n \geq 0$ we 
define a complete  flag $F_n = (V^r_n)_{r \geq 0}$ as follows: 
\begin{itemize} 
\item   $V_n^r =\epsilon^r$ for  $r \leq n$, 
\item    $V_n^{n+2p} =\epsilon^n\oplus(\eta\oplus \epsilon)^p$  for  $ p \geq 0$, 
\item   $ V_n^{n+2p+1} =\epsilon^n\oplus(\eta\oplus \epsilon)^p \oplus\eta$   for $ p  \geq 0$.
\end{itemize} 
In many cases we will work with the so-called {\em alternating flag}  $F_1$, whose subquotients \(V_1^{r+1}/V_1^{r}\), \(r\geq 0\),   alternate between \(\epsilon\) and \(\eta\), starting with $\epsilon$. 
It has been studied before in \cite[Appendix C]{CGK00}.

We denote by \(d_i\in k\) the coefficients of the expansion of \(y(\epsilon)\) in the topological basis induced by $F_0$,
\begin{equation*}
  y(\epsilon)=\sum_{i\geq 0} d_i \cdot y(V_0^i).
\end{equation*}
By applying the action $l_{\eta}$ to both sides we obtain
\begin{equation*}
  y(\eta)=\sum_{i\geq 0} d_i \cdot y(V_1^i).
\end{equation*}
Note that $d_0 = e$ by definition of the Euler class. 

We will now study the coordinate change induced by passing from the flag  $F_{n+1}$  to the flag $F_n$. For $n \geq 0$ a basis element induced by $F_{n+1}$ has one of the following forms: 
\begin{itemize} 
 \item $y(\epsilon), \dots,y(\epsilon)^n$,
 \item $ y(\epsilon)^{n+1}(y(\eta)y(\epsilon))^{p} \text{ with } p \geq 0$, 
 \item $ y(\epsilon)^{n+1}(y(\eta)y(\epsilon))^p y(\eta) \text{ with } p\geq 0$. 
\end{itemize} 
The basis elements of the first and last type are also part of the basis induced by $F_n$. 
Therefore we only have to express the basis elements of the second type in the basis induced by $F_n$,
\begin{align*}
  y(\epsilon)^{n+1}(y(\eta)y(\epsilon))^{p}&=y(\epsilon)^{n}(y(\eta)y(\epsilon))^{p}\sum_{i}d_i \cdot y(V_0^i)\\
&=y(\epsilon)^{n}(y(\eta)y(\epsilon))^{p}\sum_{i}\big( d_{2i} \cdot (y(\epsilon)y(\eta))^i+ d_{2i+1} \cdot (y(\epsilon)y(\eta))^iy(\eta)\big)\\
&=\sum_{i}\big( d_{2i} \cdot y(\epsilon)^n (y(\epsilon)y(\eta))^{i+p}+ d_{2i+1} \cdot y(\epsilon)^n (y(\epsilon)y(\eta))^{i+p}y(\eta) \big).
\end{align*}
This implies the following coordinate change formula. 

\begin{lemma} \label{lem:base_change} Let $n \geq 0$ and 
\begin{equation*} \label{eq:develop}
      \sum_{i \geq 0}  \gamma^{n+1}_i \cdot y(V_{n+1}^i) = \sum_{i \geq 0}  \gamma^{n}_{i}  \cdot y (V_{n}^i)  \in R, 
\end{equation*} 
with coefficients  $\gamma^{n+1}_i$ and $\gamma^n_i$  in $k$. Then we have \begin{itemize} 
   \item $\gamma^n_i = \gamma^{n+1}_i $ for $i < n$,
   \item $  \gamma^n_n = \gamma^{n+1}_n +  e \gamma^{n+1}_{n+1}$.  
\end{itemize} 
\end{lemma} 

Given two flags $F_n$  and $F_m$, where $n,m \geq 0$, we have a topological basis $(y(V_{n}^i)\otimes y(V_{m}^j))_{i,j \geq 0}$ of the complete $k$-module \(R\widehat{\otimes} R\).
We denote by \(\beta _{i,j}^{n,m} \in k \) the coefficients in the expansion of \(\Delta(y(\epsilon))\) with respect to this basis, 
\begin{equation} \label{introbeta} 
  \Delta(y(\epsilon))=\sum_{i,j \geq 0} \beta_{i,j}^{n,m} \cdot y(V_{n}^i)\otimes y(V_{m}^j).
\end{equation}
Note that we can study  coordinate change formulas seperately for the index pairs $(n,i)$ and $(m,j)$ while fixing the other 
index pair. 
Keeping Lemma \ref{lem:base_change} in mind we will now introduce some special elements in $k$:
\begin{itemize}
\item For \(0 \leq i<n\) and \(0 \leq j<m \) the coefficient \(\beta_{i,j}^{n,m}\) is independent of \(n, m\) and we denote this element by \(\alpha_{i,j}\).
\item For $0 \leq n$ and \(0 \leq  j<m\) the coefficient \(\beta_{n,j}^{n,m}\) is independent of \(m\) and we denote this element by \(\sigma_{n,j}\).
\item For $0 \leq m$ we set  \(\tau_{m} =\beta_{0,{m}}^{0,m}\). 
\end{itemize}

\begin{lemma} \label{lem:relations} The elements $\alpha_{i,j}$, $\sigma_{n,j}$ and $\tau_m$ satisfy the following relations:
  \begin{itemize} 
     \item[a)] $\tau_0 = 0$, \(\sigma_{0,0}=e\), \(\sigma_{1,0}=1\) and \(\sigma_{n,0}=0\) for \(n> 1\), 
     \item[b)] $\sigma_{n,j} - \alpha_{n,j} = e \sigma_{n+1, j}$ for all $j, n \geq 0$,  
     \item[c)] $ \tau_{m}  - \sigma_{0,m} =  e\tau_{m+1}$ for all $m \geq 0$. 
   \end{itemize} 
\end{lemma}
\begin{proof} By the equivariance of the comultiplication with respect to the action $l$, and using $l_{\eta}^2 = l_{\eta^2} = \id$, we have 
\[
     \Delta(y(\epsilon))=\sum_{i,j \geq 0} \beta_{i,j}^{1,1} \cdot l_{\eta} y(V_{1}^i)\otimes l_{\eta} y(V_{1}^j) = \sum_{i,j \geq 0} \beta_{i,j}^{1,1} \cdot y(V_{0}^i)\otimes y(V_{0}^j) \, . 
\]
We conclude $\tau_0 = \beta_{0,0}^{0,0} = \beta_{0,0}^{1,1} = 0$, the last equation by the co-unitality of the coproduct. 
Also, for \(n\geq 0\), we have
  \begin{align*}
    y(\epsilon)&=l_\epsilon(y(\epsilon))=(\theta(\epsilon)\otimes \Id)\circ \Delta(y(\epsilon))\\
    &=(\theta(\epsilon)\otimes \Id)\big(\sum_{i,j} \beta^{n+1,n}_{i,j} y(V_{n+1}^i)\otimes y(V_n^j)\big)\\
    &= \sum_j \beta^{n+1,n}_{0,j} y(V_n^j).
  \end{align*}
Moreover we have $\sigma_{n,0}=\beta^{n,n+1}_{n,0} = \beta^{n+1,n}_{0,n}$, the first equation by definition, the second 
equation by symmetry of the coproduct. 
Since for \(n\geq 1\) we have \(y(V_n^1)=y(\epsilon)\), and \(y(\epsilon)=\sum_i d_i \cdot y(V_0^i)\) with \(d_0=e\) by the definition of the Euler class, the remaining parts of assertion a) follow. 

Assertion b) and c) follow from the second coordinate change formula in Lemma \ref{lem:base_change}, which implies, for $n \geq 0$ and $j < m$, 
\[
   \alpha_{n,j} =  \beta^{n+1,m}_{n,j} = \beta^{n,m}_{n,j} - e \beta^{n+1, m}_{n+1,j} = \sigma_{n,j} - e\sigma_{n+1, j}  \, , 
\]
and for all $m  \geq 0$
\[
   \sigma_{0,m} = \beta^{0,m+1}_{0,m} = \beta^{0,m}_{0,m} - e\beta^{0,m+1}_{0,m+1} = \tau_{m} - e \tau_{m+1}. 
\]

\end{proof}

\begin{definition} \label{def:tame} The given equivariant formal group law  is called {\em tame} if $d_i = 0$ for all $i > 1$ and \(d_1=1\), i.e. if 
\[
    y(\eta) = y(\epsilon) + e \textrm{ and } y(\epsilon) = y(\eta) + e . 
\]
\end{definition} 
Note that for tame equivariant formal group laws we have $2 e = 0$. The additive equivariant formal 
group law (cf.~\cite[Appendix A]{CGK00}) is tame. Another tame group law will be described in Example \ref{exotic} below.

\begin{lemma} \label{lem:tame_base_change} For a tame equivariant formal group law the coefficients appearing in 
Lemma \ref{lem:base_change} 
satisfy the following relations: 
\begin{equation*} 
  \gamma^{n+1}_{i} =  \begin{cases}  \gamma^{n}_{i} \textrm{, if } i < n  \textrm{ or } i \not\equiv n \mod 2 \, ,  \\ 
                                                                 \gamma^{n}_{i} - e \cdot \gamma^{n}_{i+1} \textrm{, if  } i \geq n \textrm{ and } 
                                                                    i \equiv n  \mod 2 \, . 
                                                                 \end{cases} 
\end{equation*} 
\end{lemma} 

\begin{proof} This follows by the calculation preceding Lemma \ref{lem:base_change} together with the 
fact that $d_0 = e$, $d_1 = 1$ and $d_i = 0$ for all $i > 1$ for tame equivariant group laws. 
\end{proof} 

From this we derive the following coordinate change formula.

\begin{lemma} \label{lem:base_change_comp} Let the given equivariant formal group law be tame. 
If $i \geq n$, then in the formula
\[
  \gamma^{n+1}_{i} = \sum_{\ell = 0}^{n } x_{i,n,\ell} \cdot e^{\ell} \cdot \gamma^1_{i+ \ell} , 
\] 
which is implied by Lemma \ref{lem:tame_base_change},  the coefficients $x_{i,n,\ell}$ satisfy the congruence 
\[
   x_{i ,n,\ell} \equiv \binom{\ell+[(n - \ell)/2]}{\ell} \mod 2 , 
\]
if $i + \ell$ is even. Here $[ - ] $ denotes the Gau\ss{} bracket. 
\end{lemma} 

\begin{proof} 
The assertion is clear for $n =0$. In the induction 
step we assume the assertion holds for $n = n_0$. Assume that $i \geq n_0 +1$, $0 \leq \ell \leq n_0 +1$ and $i+ \ell$ is even. We distinguish the following cases: 
\begin{itemize} 
   \item $n_0$ and $i$ (and hence also $n_0$ and $\ell$) have the same parity. 
   \item $n_0$ and $i$ (and hence also $n_0$ and $\ell$) have different parities. 
\end{itemize} 
In the first case  we have 
\[
      \gamma^{n_0+2}_i = \gamma^{n_0 +1}_{i} 
\]
such that by the induction assumption and the fact that $n_0$ and $\ell$ have the same parity 
\[
     x_{i, n_0+1, \ell} = x_{i,n_0, \ell} \equiv  \binom{\ell+[(n_0 - \ell)/2]}{\ell} = \binom{\ell+[((n_0+1)- \ell)/2]}{\ell} \mod 2,
\]
completing the induction step. In the second case we first notice that the assertion of the lemma is 
clear for $\ell = 0$. In the case $\ell > 0$ we have
\[
      \gamma^{n_0+2}_i = \gamma^{n_0 +1}_{i} - e \gamma^{n_0+1}_{i+1} \equiv \gamma^{n_0+1}_i + e \gamma^{n_0+1}_{i+1} 
      \mod 2
\]
such that, again by the induction assumption, 
\[
    x_{i, n_0+1, \ell} = x_{i, n_0, \ell} + x_{i+1, n_0, \ell-1} \equiv \binom{\ell+[(n_0 - \ell)/2]}{\ell} + \binom{\ell-1+[(n_0 - (\ell-1))/2]}{\ell-1} \mod 2. 
\]
Using the fact that $n_0$ and $\ell$ have different parities  the last sum is equal to 
\[
\binom{\ell-1 +[((n_0 +1) - \ell)/2]}{\ell} + \binom{\ell -1 +[((n_0+1) - \ell)/2]}{\ell-1} = \binom{\ell +[((n_0+1) - \ell)/2]}{\ell}, 
\]
completing the induction step in the second case.

\end{proof}

\section{A section of the classifying map $\lambda_{\Z/2} : L_{\Z/2} \rightarrow \MU^{\Z/2}_*$} \label{section} 

Let $A = \Z/2$ and let the coproduct of the universal non-equivariant  formal 
group law be given by $\Delta(z) = \sum_{ij} a_{ij} \cdot z^i \otimes z^j$, where the elements 
$a_{ij}$, $i, j \geq 0$, generate the non-equivariant Lazard ring $L$ \cite{Laz55}. 
By \cite[Section 2]{Str01} the coefficient ring $\MU^{A}_*$ of $A$-equivariant unitary bordism 
is given as an algebra over $L$ by 
generators $s_{nj}$, $n, j \geq 0$, and $t_m$, $m \geq 0$, and relations 
 \begin{itemize} 
    \item $t_0 = 0$, $s_{10} = 1$ and  $s_{n0} = 0$ for $n > 1$, 
    \item $s_{nj} - a_{nj} = e s_{n+1, j}$, 
    \item $t_m - s_{0m} = et_{m+1}$.
 \end{itemize} 
 Here $e$ is an abbreviation for $s_{00}$, and this element corresponds to the Euler class in $\MU^{A}_{-2}$ associated to the representation $\eta$. 
 
 \begin{example} \label{exotic} Introducing the additional relations 
 \begin{itemize} 
    \item $a_{ij} = 0$ for $ i + j \geq 2$,
    \item  $s_{01} = 1$, $s_{0j}=0$ for $j \geq 2$, $s_{nj} = 0$ for $j \neq 2$ and $n \geq 1$, 
     \item $t_1 = 1$, $t_{m}= 0$ for $m \geq 1$, 
 \end{itemize} 
 we obtain a tame $\Z/2$-equivariant formal group law with a representing 
 ring which is given as a $\Z[e]/(2e)$-algebra by generators
 $s_{n2}$, $n \geq 1$, and relations $e s_{12} = 0$, $s_{n 2} = e s_{n+1, 2}$ for 
 $n \geq 1$. When viewed as a graded equivariant formal group law all elements of positive 
 degree in this ring are infinitely $e$-divisible and $e$-torsion. 
 \end{example} 

We now combine this description of $\MU^{A}_*$ with the calculus developed in Section 
\ref{sec:base_change}.

\begin{prop} The assignment 
\begin{itemize} 
  \item $a_{ij} \mapsto \alpha_{ij}$, 
  \item $s_{nj} \mapsto  \sigma_{nj}$, 
  \item  $t_{m}  \mapsto \tau_{m}$ 
 \end{itemize} 
defines a graded ring map $\mu_{A} : \MU^A _* \to L_A$ which satisfies $\lambda_A  \circ \mu_A = \id$. 

\end{prop} 

\begin{proof}
By the relations for the generators for  $\MU_*^A$ given by Strickland and by Lemma \ref{lem:relations} 
above we indeed get a well defined ring map $\MU^A_* \to L_{A}$. 

It remains to check that the canonical map $\lambda_A : L_{A}  \to \MU^A_*$ sends the elements $\alpha_{ij}$, 
$\sigma_{nj}$ and $\tau_{m}$ to the elements $a_{ij}$, $s_{nj}$ and $t_{m}$. Consider the commutative diagram 
\[
  \xymatrix{   L_{A}  \ar[r]^{\lambda_A} \ar[d] & \MU^A_* \ar[d] \\
                       \widehat{L_{A}} \ar[r]^{\widehat{\lambda_A}} & \widehat{\MU^A_*}    }
 \]
relating the map $\lambda_A$ to the induced map of  completions at the ideal $(e)$. 

The map $\widehat{\lambda_A}$ can be identified with the identity $L[[e]]/[2](e) \to L[[e]]/[2](e)$, using the canonical  
isomorphism $\widehat{L_{A}} \cong L[[e]]/[2](e)$ 
from \cite[Cor. 6.6]{G01}, which is induced by  $\alpha_{ij} \mapsto a_{ij}$, $ e \mapsto e$, 
and the $L$-algebra 
isomorphism $\widehat{\MU^A_*} \cong L[[e]]/[2](e)$ from \cite[Section 4]{Str01}. 

Furthermore the completion map 
\[
   \MU^A_* \to \widehat{\MU^A_*} \cong L[[e]]/[2](e) = \MU^{-*}(B \Z/2) 
\]
appearing as the right hand vertical map in the above diagram can be identified  
with a ``bundling map'' of tom Dieck, which, in the case relevant for us, was shown to be injective 
in \cite[Prop. 6.1 and preceding explanations]{tD70}. Alternatively the injectivity of the completion map follows from the 
results in \cite{Str01}.

 We therefore need to show that $\sigma_{nj}$ and $s_{nj}$ on the one hand, 
 and $\tau_{m}$ and $t_{m}$ on the other, are mapped to the same elements under the left and right hand vertical 
 maps in the above diagram. By the recursive formulas 
for $\sigma_{nj}$ and $\tau_{m}$ from Lemma \ref{lem:relations} and the corresponding formulas 
for $s_{nj}$ and $t_m$ from \cite[Section  4]{Str01} we arrive at the equation 
\[
     \sigma_{nj} = \sum_{\ell \geq 0} a_{n+\ell, j} e^{\ell} = s_{nj} \in L[[e]]/[2](e)
\]
and this implies, in a similar way,
\[
    \tau_{m}  = \sum_{\ell \geq 0} \sigma_{0,m+\ell} e^{\ell} = \sum_{\ell \geq 0} s_{0,m+\ell} e^{\ell} = t_{m} \in L[[e]]/[2](e) . 
\]
\end{proof}

\section{Proof of the main theorem} \label{strategy} 

In this section we explain  the proof of  Theorem \ref{main}. 
 Let us write the coproduct of the 
equivariant formal group law of $\Z/2$-equivariant unitary bordism as 
\[
    \Delta(y(\epsilon)) = \sum_{i,j\geq 0} \beta_{ij} \cdot  y(V^i) \otimes y(V^j) , 
\]
with $\beta_{ij} \in \MU^{\Z/2}_*$, where we use the alternating flag $(V^r) = (V_1^r)$. Setting 
\[
     \gamma_{ij} := \mu_{\Z/2}( \beta_{ij})  \in L_{\Z/2} ,
\]
with the map  $\mu_{\Z/2} : \MU^{\Z/2}_* \to L_{\Z/2}$ from the previous section, 
the coproduct of the universal $\Z/2$-equivariant formal group law takes the form
\[
    \Delta(y(\epsilon)) = \sum_{i,j}  ( \gamma_{ij} + \rho_{ij}) \cdot y(V^i) \otimes y(V^j) \, . 
\]
This defines new structure constants $\rho_{ij} \in L_{\Z/2}$ measuring the deviation 
from $\mu_{\Z/2}$ being surjective (and hence from $\lambda_{\Z/2}$ being injective). Note that 
$\rho_{ij} = 0$ for $i+j \leq 1$ and $\rho_{ij} \in \ker \lambda_{\Z/2}$ for all $i,j$. 
Hence each  $\rho_{ij}$ is infinitely $e$-divisible and \(e\)-torsion \cite{G01}.
In particular, 
\[
   \rho_{ij} \cdot \rho_{pq} = 0 
\]
for all $i,j,p,q \geq 0$. 

\begin{lemma} \label{kerlam} The kernel of the canonical map $\lambda_{\Z/2} : L_{\Z/2} \to \MU^{\Z/2}_*$ is equal to the
square zero ideal generated by $\rho_{pq}$, $p,q \geq 0$. 
\end{lemma} 

\begin{proof} It is clear that the given ideal is contained in $\ker \lambda_{\Z/2}$. Conversely, note 
that $L_{\Z/2}$ is generated as an $\Z[e, \gamma_{ij}]$-module by the elements $1$ and $\rho_{pq}$. 
If $x \in L_{\Z/2}$ lies in the kernel of $\lambda_{\Z/2}$, then the coefficient of $1$  in some expansion 
of $x$ as a linear combination of these generators 
is equal to $0$, because each $\rho_{pq}$ lies in $\ker \lambda_{\Z/2}$ and $\lambda_{\Z/2}$
in injective on $\Z[e, \gamma_{ij}] \cdot 1 \subset L_{\Z/2}$. 
Hence $x$ lies in the ideal generated by 
the elements $\rho_{pq}$. 
\end{proof}

The proof of Theorem \ref{main} is based on the following two results, the first of which 
provides an efficient  estimate of the order of the $e$-power torsion of $\rho_{ij}$. 

\begin{theorem} \label{step1} We have 
\[
    e^{i+j+ 1} \rho_{ij} = 0 
\]
for all $i, j \geq 0$. 
\end{theorem} 

This is proven in Section \ref{Eulerone}, where we introduce and investigate a  {\em normalization functor}, 
which turns  any $\Z/2$-equivarant formal group law into a $\Z/2$-equivariant formal 
group law with Euler class equal to $1$. 

Let $a_{ij}\in L$ denote the structure constants of the universal non-equivariant  formal group law, considered 
as elements in $\MU^{\Z/2}_*$ as in Section \ref{section}. 
Recal $a_{ij} = a_{ji}$, $a_{0j} = \delta_{1j}$ and $a_{ij}$ carries a grading equal to $2(i+j-1)$. 
Next, let $J \subset MU^{\Z/2}_*$ be the ideal generated by $a_{ij}$, $s_{nj}$, and $t_m$, with 
$i+j \geq 2$, $n + j  \geq 2$, and $m \geq 2$. In particular this ideal is generated by elements in 
strictly positive degrees. 
By the calculations in \cite{Str01} the remaining generators of $\MU^{\Z/2}_*$, as an algebra over $L[e]$, 
satisfy the relations $t_1 = 1 + e(s_{11} +t_2)$ and $s_{01} = t_1 - et_2 = 1 + es_{11}$, and hence we get    
\[
    \MU^{\Z/2}_*  / J \cong \Z[e] / (2e) . 
\]
Let $\overline{J} \subset L_{\Z/2}$ be the ideal generated by $\mu_{\Z/2}(J) \subset L_{\Z/2}$, or, in 
other words, generated by $\alpha_{ij}$, $\sigma_{nj}$, and $\tau_m$, with $i+j \geq 2$, $n + j \geq 2$, and 
$m \geq 2$. We obtain an induced $\Z/2$-equivariant formal group law with representing ring 
$L_{\Z/2} / \overline{J}$.
In Section \ref{new} we will 
show by an explicit computation that this  is in fact the additive $\Z/2$-equivariant formal 
group law. This implies the following fact. 

\begin{theorem} \label{step2} $L_{\Z/2} / \overline{J} \cong \Z[e]/(2e)$. 
\end{theorem} 

After these preparations we are in a position to prove Theorem \ref{main}.
Let some $i,j \geq 0$ be given. We claim $\rho_{ij} = 0$. This holds 
for $i+j \leq 1$ by the co-unitality of the coproduct on $L_{\Z/2}$. 
We therefore can assume $i+j \geq 2$, which implies  that the 
degree of $\rho_{ij}$ is positive. 

We abbreviate $\lambda_{\Z/2}$ by $\lambda$ and $\mu_{\Z/2}$ by $\mu$. 
Since  the degree of $\rho_{ij}$ is positive, Theorem \ref{step2} implies
\[
   \rho_{ij} = \sum_{\ell}  \gamma_{\ell} \cdot x_{\ell}  
\]
with a finite sum on the right hand side, where each $\gamma_ {\ell} \in \mu(J)$ and $x_{\ell} \in L_{\Z/2}$. Let us define 
\[
    \delta_{\ell} := \mu  \circ \lambda (x_{\ell}) \textrm{ and } x_{\ell}':= x_{\ell} - \delta_{\ell} \, . 
\]
Because each $x_{\ell}' \in \ker \lambda $ we have $ \lambda  (\sum  \gamma_{\ell} \cdot x_{\ell}') = 0$ 
and hence 
\[
   \sum  \gamma_{\ell}  \cdot \delta_{\ell} =  ( \mu \circ \lambda)  ( \sum \gamma_{\ell} \cdot \delta_{\ell}) = 
  ( \mu \circ \lambda)  ( \sum  \gamma_{\ell} \cdot x_{\ell}) = \mu \circ  \lambda (\rho_{ij}) = 0. 
\]
The first equation holds because $\sum \gamma_{\ell} \cdot \delta_{\ell} \in \textrm{ im } \mu$. We conclude 
\[
   \rho_{ij} = \sum_{\ell}  \gamma_{\ell} \cdot x_{\ell}' 
\]
where each $x_{\ell}'$ is in the ideal $\ker \lambda$, which is equal to  the ideal generated by the 
elements $\rho_{pq}$ by Lemma \ref{kerlam}. Repeating this process several times we conclude 
that  for each $N > 0$ there is a relation of the form 
\[
   \rho_{ij}  = \sum_{p,q}  c_{pq} \cdot \gamma_{pq} \cdot \rho_{pq} 
\]
where $c_{pq} \in L_{\Z/2}$ and $\gamma_{pq} \in \mu(J)^{N}$ the $N$-th power of $\mu(J)$. 
Because each generator of $J$ has degree at least $2$ we conclude, 
by comparing degrees of the right and left hand sides of the last equation, that $c_{pq} \cdot \gamma_{pq}$ 
must be divisible by $e^{p+q-1 + N - (i+j-1)}$.  For $N =  i+j+1$ the exponent satisfies 
\[
   p+q-1 + N - (i+j-1) =  p+q+1 \, .  
\]
Hence for $N = i + j + 1$ we get $c_{pq} \cdot \gamma_{pq} \cdot \rho_{pq} =0$ for all $p,q$ by Theorem \ref{step1}. We must 
therefore have $\rho_{ij} = 0$, as required.

\section{Normalization functor} \label{Eulerone} 

Let $(k, R, \Delta, \theta, y(\epsilon))$ be a (graded or ungraded) 
$\Z/2$-equivariant formal group law. As before the Euler class 
is denoted by $e$, which sits in degree $-2$, if $(k,R)$ happens to be graded. 
Passing to the 
quotient ring $k / (e-1)$ we obtain a new, ungraded $\Z/2$-equivariant formal group law with Euler class equal to one. 
In this section we will present a different way to associate to $(k, R, \Delta, \theta, y(\epsilon))$ a  $\Z/2$-equivariant formal group 
law $(k', R', \Delta' , \theta' , y'(\epsilon))$  with the following properties: 
\begin{itemize} 
   \item $k'$ is a subring of $k/\Ann(e^2)$, where $\Ann(e^2)$  is the annihilator ideal of the multiplication with $e^2$. 
    If $(k,R)$ is graded, then $k'$ is concentrated in degree $0$. 
   \item The Euler class of $(k', R')$ is equal to $1$. 
   \item The construction is functorial in $(k,R, \Delta, \theta, y(\epsilon))$. 
   \item For the formal group law $( k = \MU^{\Z/2}_*, R)$, associated to  
   $\Z/2$-equivariant unitary bordism, 
   the formal group law $(k', R')$ is the universal $\Z/2$-equivariant formal group law with Euler class equal to $1$. 
\end{itemize} 

\begin{definition} We call $(k', R', \Delta' , \theta' , y'(\epsilon))$  the {\em normalization} of $(k, R, \Delta, \theta, y(\epsilon))$. 
\end{definition}

Our construction is based on the description of equivariant formal group laws relative to flags, see \cite[Section 12]{CGK00}. We work with the alternating flag $(V^r)  = (V_1^r)$ throughout. 
Let the coproduct of $(k, R)$ be given by 
\[
    \Delta(y(V^r)) = \sum_{i,j \geq 0} f^{(r)}_{i,j} \cdot y(V^i) \otimes y(V^j). 
\]
Furthermore let 
\[
    y(\eta) = \sum_{i \geq 0} d_i y(V^i)  = e + \sum_{i \geq 1} d_{i}  y(V^i) . 
\]
Now consider the free topological $k$-module $T$ with topological basis 
$(z(V^i))_{i \geq 0}$ where $z(V^0) := 1$ and define the elements 
\[
  z(\epsilon) := z(V^1)  \textrm{ and }  z(\eta) := 1 + \sum_{i \geq 1} e^{i-1}  d_{i}   z(V^{i} )   
\]
of $T$. We define a $k$-bilinear multiplication on $T$ as follows. If at least one of $r$ and $s$ is even, then we set 
\[
     z(V^r) \cdot z(V^s) : = z(V^{r+s}) . 
\]
Furthermore, we set 
\[
   z(V^1) \cdot z(V^1) := z(V^1) +  \sum_{i \geq 1} e^{i-1} d_i  z(V^{1+i})
\]
and define 
\[
    z(V^{2p+1}) \cdot z(V^{2q+1}) := z(V^{2p+2q}) z(V^1) z(V^1) . 
 \]
This determines a continuous, $k$-bilinear product on all of $T$. It is 
 easy  to check that it is commutative and associative, where the 
 last point follows from the equation
 \[
     \big( z(V^1) \cdot z(V^1)\big) \cdot z(V^1) = z(V^1) \cdot \big( z(V^1) \cdot z(V^1) \big) 
 \]
 which is implied by commutativity of the product. 
 
 We first have to check that  the topological ring $T$ satisfies the conditions (Flag) and (Ideal) 
 from \cite[Section 12]{CGK00}. By definition we have 
\[
    z(V^{2p}) z(\epsilon) = z(V^{2p}) z(V^1) = z(V^{2p+1}) 
\]
which is consistent with (Flag), whereas by definition of $z(\eta)$ 
\[
   z(V^1) z(\eta) = z(V^1) ( 1 + \sum_{i \geq 1} e^{i-1} d_{i}  z(V^{i} ) )  . 
\]
If this element is equal to $z(V^2)$, then the condition (Flag) is satisfied on $T$ by the definition of the product on $T$. 

For this and later calculations we set
\[
   \overline{z}(V^r) : = e^r z(V^r) \in T
\]
for all $r \geq 0$ and consider the 
map $\Psi : R \to T$ of topological $k$-modules  given on basis elements by 
\[
    y(V^r) \mapsto \overline{z}(V^r) . 
\] 
If  the original formal group law $(k,R)$ is 
graded and we consider all  $z(V^r) \in T$ as sitting in degree $0$, 
then the degree of $\overline{z}(V^r)$ is equal to $-2r$ and the map 
$\Psi$ is grading preserving. Moreover  the definition of the ring structure on 
$T$ together with the calculation 
\[
  \Psi(y(V^1) \cdot y(V^1)) = \Psi \big( e y (V^1) + \sum_{i \geq 1} d_i  y (V^{1 + i})  \big)= 
   e \overline{z}(V^1) + \sum_{i \geq 1} d_i \overline{z}(V^{1 + i}) = \overline{z}(V^1) \cdot \overline{z}(V^1) 
\]
shows that the map $\Psi$ is a ring map. Note that the elements 
$\overline{z}(\epsilon) := e z(\epsilon)$ and 
\[
    \overline{z}(\eta) := e z(\eta) = \sum_{i \geq 0} d_i \overline{z} (V^i ) 
\]
are the images of $y(\epsilon)$ and $y(\eta)$ under the map $\Psi$. 
After these preparations we can calculate
\[
     e^2 \cdot z(V^1) z(\eta) = \overline{z}(V^1) \overline{z}(\eta) = \Psi ( y(V^1)) \cdot \Psi (y(\eta)) = 
     \Psi( y(V^1) \cdot y(\eta) ) = \Psi ( y(V^2)) = 
     e^2 \cdot z(V^2) . 
\]
The fourth equation follows from the relation (Flag) in $R$. 
In other words: The coefficients of  $z(V^1) z(\eta) - z(V^2) \in T$ are in the annihilator ideal $\Ann(e^2) \subset k$. 
Hence the normalization condition (Flag) is satisfied in $T$, if we pass 
from the ring $k$ to the quotient ring $k/\Ann(e^2)$. In this case  also the condition (Ideal) 
follows immediately. 

Next we define an $A^*$-action on the topological ring $T$. We set $l_{\epsilon} := \id$, 
\[
    l_{\eta} z(V^{2p} ) := z(V^{2p})\, , 
\]
 and 
\[
   l_{\eta} z(V^{2p+1})  := z(V^{2p}) \cdot z(\eta) . 
\]
 This map is extended $k$-linearly onto $T$. We need to examine the compatibility of 
 $l_{\eta}$ with the multiplication on $T$ defined before as well as the 
 property $l_{\eta}^2 = \id$. 

By definition we have 
\[
    l_{\eta} ( z(V^r) \cdot z(V^s)) = l_{\eta} z(V^r) \cdot l_{\eta} z(V^s) \, , 
\]
if either $r$ or $s$ is even. If both $r = 2p+1$ and $s = 2q+1$ are odd, then, by definition,
\[
  l_{\eta} (z(V^r) \cdot z(V^s)) = l_{\eta}\big( z(V^{2p+2q} )  \cdot z(V^1) \cdot z(V^1) \big) 
  = z(V^{2p+2q}) \cdot l_{\eta} \big( z(V^1) \cdot z(V^1) \big)
\]
and 
\[
    l_{\eta} (z(V^r) ) \cdot l_{\eta} (z(V^s))  = z(V^{2p+2q}) \cdot l_{\eta} z(V^1) \cdot l_{\eta} z(V^1) . 
\]
Thus we need  to show 
\[
   l_{\eta} \big( z(V^1) \cdot z(V^1)) =  l_{\eta} ( z(V^1)) \cdot l_{\eta} (z(V^1))  . 
\]
For this we calculate 
\[
   e^2 \cdot l_{\eta} \big( z(V^1) \cdot z(V^1)) = 
   l_{\eta} \big( \overline{z} (V^1) \cdot \overline{z}(V^1)) = 
   l_{\eta}( \overline{z}(V^1)) \cdot l_{\eta} (\overline{z}(V^1)) = 
   \overline{z}(\eta) \cdot \overline{z}(\eta) = 
   e^2 \cdot z(\eta) z(\eta) . 
\]
The second equation uses the fact that the above ring map $\Psi$ is compatible with 
the map $l_{\eta}$. For this assertion we notice that indeed 
$l_{\eta}(\overline{z}(\epsilon)) = e l_{\eta} (z(\epsilon)) = e z(\eta) = \overline{z}(\eta)$ and 
\[
   l_{\eta} ( \overline{z}(\eta) ) = l_{\eta} ( \sum_{i \geq 0} d_i \overline{z}(V^i)) = e + \overline{z}(\eta) ( 
    \sum_{i \geq 1} d_i \overline{z}(V^{i-1}) ) = \overline{z}(\epsilon), 
\]
where the last equation follows from the corresponding relation in $R$ and application of the ring map $\Psi$. In summary we see that $l_{\eta}$ is 
a ring map on $T$ after passing to the  coefficient ring $k/\Ann(e^2)$. 

The equation $l_{\eta} \circ l_{\eta} = \id$ on $T$ holds after passing to the coefficient ring $k/\Ann(e^2)$, because 
\[
  e  \cdot l_{\eta} z(\eta) =  l_{\eta} \overline{z}(\eta) = \overline{z}(\epsilon) = e \cdot z(\epsilon) 
 \]
 and $l_{\eta} z(\epsilon) = z(\eta)$, by definition. 
 
Let us  turn to the definition of the coproduct on $T$. 

\begin{lemma} \label{coeff_copr} For the structure constants of the coproduct in $R$ we have 
\[
   f^{(2)}_{0,j} = \delta_{2,j} \textrm{ and } e f^{(2)}_{1,j} = 0 
\]
for all $j \geq 0$. 
\end{lemma} 

\begin{proof}  Recall 
\[
    \Delta(y(V^2)) = \sum_{i,j \geq 0} f^{(2)}_{i,j} \cdot y(V^i) \otimes y(V^j) . 
\]
The first equality of the lemma holds by the co-unitality of $\Delta$. For the second equality we use the fact 
that the coproduct $\Delta$ is compatible with $l_{\eta}$: 
\[
   \Delta(y(V^2)) =  \Delta(l_{\eta} y(V^2)) = \sum_{i,j \geq 0} f^{(2)}_{i,j} \cdot ( l_{\eta} y(V^i) ) \otimes y(V^j) \, . 
\]
Hence the second equation in the lemma follows from $l_{\eta} (V^1) = e + \sum_{i \geq 1} d_i y(V^i)$ and 
$l_{\eta}y(V^r) \in (y(V^1))$ (the ideal spanned by $y(V^1))$ for all $r \geq 2$, by comparing coefficients and 
using the first part of the lemma. 
\end{proof} 

Now we set 
\begin{eqnarray*} 
    \Delta (z(\epsilon)) & := & \sum_{i, j \geq 0} e^{i+j-1} f^{(1)}_{i,j}  \cdot z(V^i) \otimes z(V^j) , \\
    \Delta (z(V^2)) & := & \sum_{i,j \geq 0} e^{i+j-2} f^{(2)}_{i,j} \cdot z(V^i) \otimes z(V^j) .
\end{eqnarray*} 
Note that on the right hand sides only non-negative exponents occur at $e$, by the co-unitality of the coproduct $\Delta$ on $R$, compare the first part of Lemma \ref{coeff_copr}. We now define, for $p \geq 1$, 
\[
   \Delta(z(V^{2p}))  := \Delta(z(V^2))^p  , 
 \]
 and 
 \[
   \Delta(z(V^{2p+1}))  := \Delta(z(V^{2p}) )\cdot \Delta(z(\epsilon)). 
 \]
 It follows from the definition of the multiplication on $T$ and from the second part of Lemma \ref{coeff_copr} 
 that this extends to a continuous map 
 \[
     \Delta : T \to T  \widehat{\otimes} T  
 \]
 after passing to the coefficient ring $k/\Ann(e^2)$. We also notice that this coproduct $\Delta$ is 
 compatible with the coproduct on $R$ and the ring map $\Psi : R \to T$. 
 By definition 
 \[
    \Delta(z(V^{r}) \cdot z(V^s) ) = \Delta( z(V^r) ) \cdot \Delta(z(V^s)) 
 \]
 if at least one of $r$ and $s$ is even. If $r$ and $s$ are odd, then we  only need to check 
 \[
    \Delta\big( z(V^1) \cdot z(V^1)\big) = \Delta(z(V^1)) \cdot \Delta(z(V^1)) . 
 \]
 Applying the same reasoning as before we know that this equation holds after multiplication with $e^2$, by 
 using the corresponding relation for the coproduct on $R$. Hence this equation holds 
 after passing to the coefficient ring $k / \Ann(e^2)$. 
 
 Next we check compatibility of $\Delta$ with the left $A^*$-action on $T$. Because this action is by ring maps, it is enough 
 to check \begin{itemize} 
 \item $ \Delta( z(\eta)) = (l_{\eta} \otimes \id) (\Delta(z(\epsilon)))  = (\id \otimes l_{\eta}) ( \Delta(z(\epsilon)))$, 
 \item $\Delta(z(\epsilon)) = (l_{\eta} \otimes \id) (\Delta(z(\eta))) = (\id \otimes l_{\eta}) (\Delta(z(\eta)))$, 
 \item $\Delta ( z(\epsilon)) = ( l_{\eta} \otimes l_{\eta}) ( \Delta(z(\epsilon))$, 
 \item $\Delta ( z(\eta)) = ( l_{\eta} \otimes l_{\eta}) ( \Delta(z(\eta))$. 
\end{itemize} 
All of these equation are true after multiplication with $e$, hence we are fine if we work over the coefficient ring $k/\Ann(e^2)$. 
Finally  we check co-associativity. Because $\Delta$ is multiplicative on $T$, if we work over the
coefficient ring $k/\Ann(e^2)$, it is enough to consider the equations
\[
  (  \Delta \otimes \id) \circ \Delta( z(\epsilon)) = (\id \otimes \Delta) \circ \Delta ( z(\epsilon))  ,
 \]
 and 
  \[
  (  \Delta \otimes \id) \circ \Delta  ( z(V^2)) = (\id \otimes \Delta) \circ \Delta ( z(V^2))  .
 \]
The first equation 
 holds after multiplication with $e$, and the second equation holds after multiplication with $e^2$. 
 Hence both equations hold after passing to the coefficient ring  $k/\Ann(e^2)$. 
 
In summary, using the discussion of \cite[Section 12]{CGK00},
we have defined a $\Z/2$-equivariant formal group law $(k/\Ann(e^2), T, \Delta, \theta, z(\epsilon))$, where 
we write $T$ instead of $T / (\Ann(e^2))$ by a slight abuse of notation. 
The augmentation $\theta : T \to (k/\Ann(e^2))^{A^*}$ is given by the constant 
term in the expansion  relative to the flag $(V_1^r)$ (resp. relative to the flag $(V_0^r)$), 
at the representation $\epsilon$ (resp. at the representation $\eta$).

By definition this formal group law 
has Euler class equal to $1$. Now we let $k' \subset k/ \Ann(e^2)$ be the subring generated by the coefficients 
$e^{i+j-1} f^{(1)}_{i,j}$, $i, j \geq 0$, of the coproduct on $T$ (regarded as elements in $k/\Ann(e^2)$) and define 
$R'$ as the free topological $k'$-module with basis $(z(V^r))_{r \geq 0}$. If $(k,R)$ is graded, 
then $k'$ is indeed concentrated in degree $0$. Regarding $R'$ 
as a subset of $T$ we note that the product, coproduct and $A^*$-action on $T$ restrict to corresponding 
structures on $R'$, compare \cite[Theorem 16.1]{CGK00}.
 Also the augmentation $\theta$ restricts to an augmentation $\theta' : R' \to (k')^{A^*}$. Setting 
$y'(\epsilon) := z(\epsilon)$ this concludes the construction of $(k', R', \Delta', \theta', y'(\epsilon))$. The 
functoriality of this construction is clear. 

The next result highlights an important example. 

\begin{prop} \label{norm_univ} 
The normalized formal group law $R'$ associated to  $\Z/2$-equivariant unitary bordism $k = \MU^{\Z/2}_*$ 
is the universal $\Z/2$-equivariant formal group law with Euler class $1$. 
\end{prop} 

\begin{proof} Set $k = \MU^{\Z/2}_*$ and let $R$ be the topological $k$-algebra of the associated 
$\Z/2$-equivariant formal group law. We work with the notation from \cite{Str01}, repeated in Section \ref{section} 
above. 
By \cite[Cor.~10]{Str01} the annihilator ideal $\Ann(e^2) \subset k$ 
is generated by $(t_1+1)$, and in fact equal to the annihilator ideal of the multiplication with $e$.  By Section \ref{section} 
we can identify the distinguished generators $a_{ij}$, $t_{m}$ and $s_{nj}$ of $k$ with certain coefficients of the coproduct 
$\Delta(y(\epsilon))$ in $R$ developped with respect to suitable flags. By our definition of the coproduct on $R'$ the ring $k'$ 
is therefore the subring of $k/(t_1 +1)$ generated by the elements 
\[
  \overline{a}_{ij} := e^{i+j-1} a_{ij} \, , \, \,  \overline{s}_{nj} := e^{n+j-1} s_{nj} \, , \, \, \overline{t}_m :=  e^{m-1} t_{m} ,
\]
with $i+j \geq 1$, $n+j \geq 1$, and $m \geq 1$. These elements are only subject to the relations 
\[
    \overline{t}_m - \overline{s}_{0m} = \overline{t}_{m+1} \, , \, \,  \overline{s}_{nj} - \overline{a}_{nj}= \overline{s}_{n+1, j} 
\]
for all $m$, $j$, and $n$. This implies that $k'$ is generated as a $\Z$-algebra by the elements $\overline{a}_{ij}$, $i + j \geq 2$,  and $\overline{s}_{0m}$, $m \geq 1$, where the generators $\overline{a}_{ij}$ satisfy the same relations as in the non-equivariant Lazard ring and $\overline{s}_{0m}$ are free polynomial generators. 
Using Strickland's calculation of $k$ we hence  conclude, on the one hand, that the quotient map 
$k \mapsto k / (e-1)$ induces an isomorphism 
\[
    k' \cong k/(e-1) . 
\]
On the other hand we observe that $L_{\Z/2} / (e-1)$ is the underlying ring of the universal $\Z/2$-equivariant formal 
group law with Euler class equal to $1$. But the classifying map 
\[
    \lambda_{\Z/2} : L_{\Z/2} \to k 
\]
is surjectivce and elements in the kernel are Euler torsion, see \cite{G01}. It hence 
induces an  isomorphism 
\[
    L_{\Z/2}  / (e-1) \cong k / (e-1) . 
\]
This finishes the proof of Proposition \ref{norm_univ}. 
 \end{proof} 

We are now in a position to prove Theorem \ref{step1}. Consider the classifying map 
\[
   \lambda : L_{\Z/2} \to k := \MU^{\Z/2}_*
\]
of the equivariant formal group law of 
$\Z/2$-equivariant unitary bordism and the section $\mu$ of this map 
constructed in Section \ref{section}. Applying the normalization functor 
we obtain induced maps $\lambda'   : L'_{\Z/2} \to k'$
and $\mu'   : k'  \to L'_{\Z/2}$ satisfying $ \lambda' \circ  \mu'= \id$, by functoriality.
Furthermore  we have 
\[
    \mu'( e^{i+j-1} \beta_{ij} ) = e^{i+j-1} \gamma_{ij}
 \]
by the definition of $\beta_{ij}$ and $\gamma_{ij}$ in Section \ref{strategy}. 
Proposition \ref{norm_univ} implies that  the equivariant formal group law for $L'_{\Z/2}$  is classified by a ring map 
\[
    \phi : k' \to L'_{\Z/2} . 
\] 
Because both $\lambda'$ and $\phi$ are classifying maps, 
 they are inverse to each other, and because  $ \lambda'\circ  \mu' = \id$ we in fact have  $\phi = \mu' $. 
 Since
\[
   \phi(e^{i+j-1} \beta_{ij}) = e^{i+j-1}( \gamma_{ij} + \rho_{ij})
\]
this implies 
\[
     e^{i+j-1} \rho_{ij} = 0 \in L'_{\Z/2} .
\]
Finally, because $L'_{\Z/2} \subset  L_{\Z/2} / \Ann(e^2)$, this implies the equation 
\[
    e^{i+j+1} \rho_{ij} = 0 
\]
in $L_{\Z/2}$, and hence the assertion of Theorem \ref{step1}.

\section{Computation of a particular $\Z/2$-equivariant formal group law} \label{new} 

Let $(L_{\Z/2}, R, \Delta, \theta, y(\epsilon))$ denote the universal $\Z/2$-equivariant formal group law. 
We consider the graded  ideal $\overline{J} \subset L_{\Z/2}$ defined in Section \ref{strategy}, 
spanned by the homogenous elements $a_{ij}$, $\sigma_{nj}$ and 
$\tau_m$  for 
$i+ j \geq 2$, $n+j \geq 2$,  and $m \geq 2$. The resulting $\Z/2$-equivariant formal group law $( L_{\Z/2} / \overline{J} , R / ( \overline{J}) , \Delta)$ 
has the form 
\[
     \Delta(y(V_1^1)) = y(V_1^1) \otimes 1 + 1 \otimes y(V_1^1) + \sum_{i,j \geq 0} \rho_{ij} \cdot y(V_1^i) \otimes y(V_1^j) 
\]
with respect to the alternating flag.  By abuse of notation we here denote by $\rho_{ij} \in L_{\Z/2} / \overline{J}$ the images of the 
structure constants $\rho_{ij}$  introduced in 
Section \ref{strategy}. In particular we have $\rho_{ij} = 0$ for $i +  j \leq 1$ and 
the $\rho_{ij}$ are all infinitely $e$-divisible and $e$-torsion, such that all products 
$\rho_{ij} \cdot \rho_{pq}$ are equal to $0$. Also recall that $\rho_{ij} = \rho_{ji}$ for all $i,j \geq 0$. 

The section $\mu_{\Z/2} : \MU^{\Z/2}_* \to L_{\Z/2}$
induces a section $\MU_*^{\Z/2} / J \to L_{\Z/2}/ \overline{J}$ of the canonical map 
$\lambda : L_{\Z/2}/ \overline{J} \to \MU_*^{\Z/2} / J$. We can hence consider  $\MU^{\Z/2}_* / J = \Z[e]/(2e)$ 
as a subring of $L_{\Z/2}/ \overline{J}$. In particular $2e = 0$ in $L_{\Z/2}/ \overline{J}$ and therefore 
\[
    2 \cdot \rho_{ij} = 0 
\]
for all $i,j  \geq 1$, by the $e$-divisibility of $\rho_{ij}$.  This will simplify the following computations considerably. 

Notice that the kernel of $\lambda: L_{\Z/2}/ \overline{J} \to \MU_*^{\Z/2} / J$ is generated by 
the structure constants $\rho_{ij}$.
In the remainder of this section we will show  that all $\rho_{ij} = 0 \in L_{\Z/2} / \overline{J}$. 
 This assertion implies Theorem \ref{step2}. 

The structure constants $\rho_{ij} \in  L_{\Z/2} / \overline{J}$ underly the following restrictions: 
\begin{itemize} 
  \item[(1)] The comultiplication $\Delta$ is co-associative.  
  \item[(2)] The elements $\sigma_{n,j} \in L_{\Z/2}$,  $n+j \geq 2$, and $\tau_m \in L_{\Z/2}$, $m \geq 2$, map to 
                 $0$ in the quotient ring $L_{\Z/2} / \overline{J}$. 
\end{itemize} 
At first we will explore Restriction (1), which  results in Proposition \ref{rel1}. Then, in Propositions \ref{dyadic}, \ref{sum} 
and \ref{taus} 
we will derive  implications from Restriction (2), making use of the coordinate change formulas from Section \ref{sec:base_change}. After 
these preparations the proof of $\rho_{ij} = 0$ will be completed at the end of this section. 
 
As a shorthand we use the notation $z^r := y(V_1^r) \in R / (\overline{J} )$, $r \geq 0$,  for the basis elements corresponding to the alternating flag. We warn the reader that in the ring $R / ( \overline{J} )$ we cannot assume the relation 
 $z^{r_1} \cdot z^{r_2} = z^{r_1 + r_2}$ for $r_1, r_2 \geq 1$. 

The $\Z/2$-equivariant additive formal group law $(k_a, R_a, \Delta_a, \theta_a, y_a(\epsilon))$ with representing ring $k_a = \Z[e] / (2e)$ defines structure constants $f^{(r)}_{ij} \in \Z[2] / (2e) $ for $r \geq 0$ by the equation 
\[ 
   \Delta_a (y_a(V_1^r)) =  \sum_{i,j \geq 0} f^{(r)}_{ij} y_a(V_1^i) \otimes y_a(V_1^j)  . 
\]
We now define $\rho^{(r)}_{ij} \in L_{\Z/2} / \overline{J}$ by the equation 
\[
   \Delta(z^r) = \sum_{i,j \geq 0} ( f_{ij}^{(r)} + \rho^{(r)}_{ij} ) \cdot  z^i \otimes z^j \, . 
\]
In particular $\rho^{(0)}_{i,j} = 0$, $\rho^{(1)}_{i,j} = \rho_{ij}$ for all $i,j$, and all $\rho^{(r)}_{ij}$ 
are in the kernel of $\lambda: L_{\Z/2}/ \overline{J} \to \MU_*^{\Z/2} / J$ and hence lie in the 
ideal spanned by $\rho_{pq}$, $p,q \geq 0$. In particular $2 \cdot \rho^{(r)}_{ij} = 0$ and $\rho^{(r)}_{ij}\cdot \rho^{(s)}_{pq} = 0 $ for all $i,j,r$ and $p,q,s$.

We obtain
\begin{align*}
  (\Delta\otimes \Id)\circ \Delta(z)&=(\Delta\otimes \Id)(\sum_{j,k}(f_{j,k}^{(1)} +\rho_{j,k}^{(1)})\cdot z^j \otimes z^k )\\
&=\sum_{l,m,j,k} (f_{l,m}^{(j)} + \rho_{l,m}^{(j)})(f_{j,k}^{(1)}+\rho_{j,k}^{(1)})\cdot z^l \otimes z^m \otimes z^k \\
&=\sum_{l,m,j,k} (f_{l,m}^{(j)} f_{j,k}^{(1)} + \rho_{l,m}^{(j)} f_{j,k}^{(1)}+f_{l,m}^{(j)}\rho_{j,k}^{(1)}) \cdot z^l \otimes z^m \otimes z^k.
\end{align*}
Here we use the vanishing of products of $\rho$'s. 
Similarly we have
\begin{align*}
  (\Id \otimes \Delta)\circ \Delta(z)&=\sum_{l,m,j,k} (f_{k,m}^{(j)} f_{j,l}^{(1)} + \rho_{k,m}^{(j)} f_{j,l}^{(1)}+ f_{k,m}^{(j)}\rho_{j,l}^{(1)} )\cdot z^l \otimes z^m \otimes z^k .
\end{align*}
Taking into account that the comultiplication $\Delta_a$ is co-associative, these calculations together with the co-associativity of $\Delta$ 
imply that  for all $k, l, m \geq 0$ we have 
\begin{equation} \label{eq:ass} 
  \sum_{\nu \geq 0}(\rho_{l,m}^{(\nu)} \cdot  f_{\nu,k}^{(1)} + f_{l,m}^{(\nu)} \cdot \rho_{\nu,k}^{(1)} )=\sum_{\nu \geq 0}( \rho_{k,m}^{(\nu)}  \cdot f_{\nu,l}^{(1)} + f_{k,m}^{(\nu)} \cdot  \rho_{\nu,l}^{(1)} ) \, . 
\end{equation}
Observe that this is a finite sum on each side, because the comultiplication $\Delta$ is continuous. 

Since the equation $2 \cdot \rho^{(r)} _{ij} = 0$ holds for all $i,j,r$, we need to compute the images of the elements $f_{pq}^{(r)}$ in $\Z/2 [e]$ 
when  evaluating Equation \eqref{eq:ass}.
In the following we use the shorthand notation $x^r := y_a(V_1^r)$, $x := x^1$. In particular $x^2 = x\cdot (x+e) = x^2 + ex$, $x^{2n} = (x^2 + ex)^n$ and 
$x^{2n+1} = x \cdot (x^2 + ex)^n$ for all $n \geq 0$. We have $ \Delta_a (x)=x \otimes 1 + 1\otimes x$ and 
\[
   \Delta_a(x^2) = \Delta_a( x) \cdot \Delta_a ( x + e) = 
   ( x\otimes 1 + 1 \otimes x ) \cdot (x  \otimes 1 + 1 \otimes x + e \cdot (1 \otimes 1)) = x^2 \otimes 1 + 1 \otimes 
   x^2 \, , 
\]
after passing to the representing ring $\Z/2 [ e]$. Hence (for even and odd $r$) we obtain
\begin{equation*}
\label{eq:binomi}
   \Delta_a(x^r) = (x \otimes 1 + 1 \otimes x)^r = \sum_{s = 0}^{r} \binom{r}{s} x^s \otimes x^{r-s} \, . 
\end{equation*}
This happens to be the same formula as for the additive non-equivariant formal group law with coordinate $x$. In $\Z/2[e]$ we hence obtain 
the equality 
\[
f_{p,q}^{(r)}=\begin{cases} \binom{p+q}{p}&\text{ if } r =  p+q,\\
0&\text{ else. }
\end{cases}
\]
Recall  that $\binom{p+q}{p}$ is equal to $0$ modulo $2$, if and only if in the binary expansions of $p$ and $q$ the digit
$1$ occurs at the same position, or, in other words, if the binary addition of $p$ and $q$ involves carryovers.

We arrive at the following conclusion resulting from Restriction (1). 

\begin{prop} \label{rel1} $L_{\Z/2} / \overline{J}$  is generated over $\Z[e]/(2e)$ by  the elements $\rho_{ij}$, 
$i+ j \geq 2$, and these elements satisfy the following relations. 
\begin{itemize} 
\item[a)]  $\rho_{ij} \cdot \rho_{p,q} = 0$ and $2 \rho_{ij} = 0$. 
\item[b)]  If $i, j \geq 1$, if either $i$ or $j$ is not a power of $2$, and if the binary addition of $i$ and $j$ involves carryovers, then
\[
   \rho_{ij} = 0  \in  L_{\Z/2} / \overline{J} . 
\]
\item[c)]   If neither $i,j\geq 1$ nor $p,q\geq 1$ fall in the case b),  and  if $i+j = p+q$, then 
\[
   \rho_{ij} = \rho_{pq} .
\]
\end{itemize} 
\end{prop} 
\begin{proof}
It remains to deal with parts b) and c). Both of them  follow from Equation \eqref{eq:ass}, where we observe that the 
elements $\rho_{l,m}^{(\nu)}$ and $\rho_{k,m}^{(\nu)}$ can only occur with a factor $1$, if $\nu = 1$, 
by our previous computation of $f^{(1)}_{p,q}$ as 
elements in $\Z/2[e]$. 

For part b) we write 
\[
    i=\sum_{\ell\geq 0} w_{i\ell} \cdot 2^{\ell} \textrm{ and }  j=\sum_{\ell\geq 0} w_{j\ell} \cdot 2^{\ell}
\]
with \(w_{i\ell},w_{j\ell}\in\{0,1\}\) where $i$, say,  is not a power of $2$. We choose $\ell_1$ with 
\(w_{i\ell_1}=w_{j\ell_1}=1\).
Then the assertion follows from Equation \eqref{eq:ass} with  \(k=j\), \(l=i-2^{\ell_1}\) and \(m=2^{\ell_1}\).

Now we turn to part c). If $i$ and $j$ are both powers of two, then under the given assumptions the same must hold for $p$ and $q$. 
We obtain  $p = i$, $p = j$ or $p = j, q =i$, and  claim c) follows  from the commutativity of the formal group law $\Delta$. 
It therefore remains to deal with the case that $i$, say,  is not a power of two. 
Let us write 
\[
   i=\sum_{\ell \geq 0} w_{i\ell} \cdot 2^{\ell}  \textrm{ and }  j=\sum_{\ell \geq 0}w_{j\ell} \cdot 2^{\ell}
\]
with \(w_{i\ell},w_{j\ell}\in\{0,1\}\) and   \(w_{i\ell}\cdot w_{j\ell}=0\) for all $\ell$. Choose $\ell_1$  with 
 \(w_{i \ell_1}=1\).  
Then it follows from Equation \eqref{eq:ass} with \(k=i-2^{\ell_1}\), \(m=2^{\ell_1}\), \(l=j\), that
\begin{equation*}
  \rho_{ij}=\rho_{i-2^{\ell_1},j+2^{\ell_1}} \, . 
\end{equation*}
In other words, we can shift the binary digit $1$ at position $\ell_1$ from the left to the right hand subscript of $\rho$. 
From this claim c) in the proposition follows. 

\end{proof}

For exploring Restriction (2) we need to work with different flags. Let us write, for $n,m \geq 1$, 
\[
     \Delta(y(V^1_1)) = y(V_n^1)  \otimes 1 + 1 \otimes y(V_m^1)   + \sum_{i,j \geq 0}\rho_{i,j}^{n,m} \cdot y(V_n^i) \otimes y(V_m^j)  \, . 
\]
Note that $y(V_1^1) = y(V_n^1) = y(V_m^1) = z^1$ by our assumption $n,m \geq 1$. 
By the co-unitality of $\Delta$ we therefore have $\rho_{i,0}^{n,m} = \rho_{0,j}^{n,m} = 0$ for all $i, j \geq 0$, 
and, using  the notation  introduced in Equation \eqref{introbeta}, we have $\rho_{i,j}^{n,m} = \beta_{i,j}^{n,m}$ 
for $i +  j \geq 2$ (notice that $\rho_{1,0}^{n,m} = \rho_{0,1}^{n,m} = 0$,  whereas $\beta_{1,0}^{n,m} = \beta_{0,1}^{n,m} = 1$ for $n,m \geq 1$). 
Also note that $\rho^{1,1}_{ij} = \rho_{ij}$ for $i,j \geq 0$,  and 
all  $\rho^{n,m}_{i,j}$ are in the 
kernel of the map $L_{\Z/2}/ \overline{J} \to \MU^{\Z/2} / J$, and hence 
lie in  the ideal generated by the elements $\rho_{pq}$.  In particular 
all   $\rho^{n,m}_{i,j}$ are $2$-torsion and arbitrary products of such 
elements vanish.

We wish to apply the   coordinate change formula in Lemma \ref{lem:tame_base_change}. 
This was originally stated for tame group laws. However, for the group law considered in this section the coefficients 
$d_i$ appearing in the base change formula preceding Lemma 
\ref{lem:base_change} are equal to those of a tame equivariant law, modulo 
elements in the kernel of $\lambda : L_{\Z/2} / \overline{J} \to \MU^{\Z/2}_*/ J$, which 
is equal to the square zero ideal generated by the elements $\rho_{pq}$. 
Hence Lemmas \ref{lem:tame_base_change}
and \ref{lem:base_change_comp} remain valid in our case of the (potentially) non-tame group law $\Delta$, if we apply it to coordinates of the form
$\gamma^{n+1}_i := \rho^{n+1, m}_{i,j}$ or $\gamma^{m+1}_j := \rho^{n,m+1}_{i,j}$, where $n,m \geq 1$ and $i,j \geq 1$. 
 Hence, for all $n, m,  j \geq 1$ we have equations
\begin{equation}
 \label{eq1} 
  \rho^{n, m}_{n,j} = \sum_{ \ell = 0}^{n-1} y_{n,\ell} \cdot e^{\ell} \cdot \rho^{1, m}_{n+ \ell, j} \, , 
\end{equation}  
where  $y_{n,\ell} \in \Z/2$, and for all $n, m , i \geq 1$ we have equations 
\begin{equation}
 \label{eq2} 
   \rho^{n,m+1}_{i,m} = \sum_{\nu = 0}^{m} x_{m, \nu} \cdot e^{\nu} \cdot \rho^{n,1}_{i, m + \nu}  ,
\end{equation}
where each $x_{m,\nu} \in \Z/2$ is equal to 
the coefficient $x_{m, m, \nu}$ from Lemma \ref{lem:base_change_comp}. Notice that 
\[
    y_{n,0} = x_{m,0} = 1 
\]
for all $n, m  \geq 1$, by the recursive formula in Lemma \ref{lem:tame_base_change}.
Let us compute the coefficients $x_{m,\nu}$ in some more cases. 

\begin{lemma} \label{lem:miracle} For all $q \geq 0$ we have 
\[
         x_{2^q, 2^q} = 1 .
\]
Furthermore, if $q \geq 2$, then for all $0 < \omega < 2^{q-1}$  we have
\[
   x_{2^{q} - \omega, \omega}  = 0 .
\]
\end{lemma} 
Note that  $x_{2^{q} - \omega, \omega}$ is the coefficient appearing in front of $e^{\omega} \cdot \rho^{n,1}_{i, 2^q}$, if 
we  develop $\rho^{n,2^q - \omega +1}_{i,2^q - \omega}$ according to Equation \eqref{eq2}. 
 The relation $x_{2^q - \omega, \omega} = 0$  will be a crucial ingredient for proving Proposition \ref{dyadic} below. 

\begin{proof}[Proof of Lemma \ref{lem:miracle}] By  Lemma \ref{lem:base_change_comp} and the discussion preceding 
Lemma \ref{lem:miracle} we have 
\[
   x_{m,\omega } = x_{m ,m,\omega} = \binom{\omega+[(m - \omega )/2]}{\omega} \mod 2 
\]
if $m + \omega$ is even. Evaluating this formula for $m = \omega = 2^q$ shows the first assertion. 
We now assume $q \geq 2$, $0 < \omega < 2^{q-1}$ and set $m = 2^q - \omega$.
Then $m+\omega$ is even because $q \geq 2$ and we obtain 
\[
   x_{m,\omega} = \binom{\omega+[( 2^q -  2 \omega)/2]}{\omega} = \binom{2^{q-1}}{\omega}  \mod 2.
\]
This vanishes because  $0 < \omega < 2^{q-1}$ by assumption. Hence the lemma 
is proven. 
\end{proof}

We can now explore Restriction (2), saying $[\sigma_{n, j}] =  [\tau_m]  = 0 \in L_{\Z/2} / \overline{J}$ for $n + j \geq 2$ and $m \geq 2$. 
Let us start with the relation $[\sigma_{n, j}] = 0$.

\begin{prop}  \label{dyadic} Let $p \geq 1$. If $1 <  j < 2^p$ is not a power of \(2\), then we have
\[
   \rho^{1,1}_{2^p, j} = 0 . 
\]
\end{prop} 

\begin{proof} 
Let \(1 <  j <2^{p-1}\) be not a power of $2$ and assume inductively that we have 
proven $\rho^{1,1}_{2^p, j'} = 0 $ for all $j < j' < 2^p$ where $j'$ is not a power of $2$ (this condition is empty for $j =2^p - 1$). 
Choose $1 \leq q \leq p$ minimal with $2^q > j$ and write $j = 2^q - \omega$ where $0 <   \omega < 2^{q-1}$.  
Using Equations \eqref{eq1} and \eqref{eq2} we obtain
\begin{equation}  \label{expansion} 
  0 =   [\sigma_{2^p, 2^{q}-\omega}] = \rho^{2^p, 2^q-\omega+1}_{2^p, 2^q-\omega} =
   \sum_{\ell=0}^{2^p -1}\sum_{\nu=0}^{2^q-\omega}  y_{2^p,\ell}\cdot x_{2^{q}-\omega,\nu}\cdot e^{\ell+\nu}\rho^{1,1}_{2^p+\ell,2^q-\omega+\nu} \, .
\end{equation} 
By part b) of Proposition \ref{rel1} $\rho^{1,1}_{2^p+\ell,2^q-\omega+\nu} \neq 0$ can only occur in one of the following cases: 
\begin{itemize} 
   \item[i)] $2^p + \ell$ and $2^{q} - \omega + \nu$ are both powers of $2$. 
   \item[ii)] The binary addition of $2^p + \ell$ and $2^{q} - \omega + \nu$ does not involve carryovers. 
\end{itemize} 
Case i) is equivalent to  $\ell = 0$ and $\nu - \omega = 0$, and the corresponding summand on the right hand side of Equation \eqref{expansion} is 
equal to $x_{2^q - \omega, \omega} \cdot e^{\omega} \cdot \rho^{1,1}_{2^p, 2^q}$ (recall that $y_{2^p, 0} = 1$). By Lemma \ref{lem:miracle} 
we have $x_{2^{q} - \omega, \omega}=0$ and hence this summand vanishes. 

Let us now assume that we are in case ii), but not in case i). We claim that $2^q - \omega + \nu + \ell < 2^p$.  In a first step we 
prove $2^q - \omega + \nu < 2^{p}$. 
Here  we notice 
$2^q - \omega + \nu < 2^{p+1}$, because $q \leq p$ and $\nu - \omega < 2^{q}$. Hence the assumption $2^q - \omega + \nu \geq 2^p$ 
together with $0 \leq \ell <  2^p$ implies that in the binary expansions of both $2^p + \ell$ and $2^q - \omega +\nu$ the digit $1$ occurs at position $p$ 
(corresponding to $2^p$), contradicting the assumption of case ii). Because $\ell \leq 2^p -1$ and the binary addition of 
$2^p + \ell$ and $2^q - \omega + \nu$ does not involve carryovers, the inequality $2^q - \omega + \nu  < 2^p$ 
 in turn implies $2^q - \omega + \nu + \ell < 2^p$, as claimed before. 

Part c) of Proposition \ref{rel1} now  implies 
\[
   \rho^{1,1}_{2^p+\ell,2^q-\omega+\nu} = \rho^{1,1}_{2^p , 2^q  - \omega + \nu + \ell} \, . 
\]
Since $2^q  - \omega + \nu + \ell$ is not a power of $2$ (by the assumption of case ii) and since we are not in case i)) 
and smaller than $2^p$ (as shown before), the last expression vanishes by our induction 
assumption, if either $\ell > 0$ or $\nu > 0$. 

In summary Equation \eqref{expansion} simplifies to $ 0 = \rho^{1,1}_{2^p, 2^q - \omega}$, finishing the induction step.

\end{proof}

\begin{prop} \label{sum} If $i, j \geq 1$ and either $i$ or $j$ is not a power of $2$, then 
\[
   \rho^{1,1}_{i, j} = 0 . 
\]
If $i$ and $j$ are both powers of $2$, then we have the relations
\[
    \rho^{1,1}_{i,j} + e^{i} \rho^{1,1}_{2i, j} = 0 \textrm{ and } \rho^{1,1}_{i,j} + e^j \rho^{1,1}_{i,2j} = 0. 
\]
\end{prop} 

\begin{proof} The first assertion follows from Proposition \ref{dyadic} and 
the parts (b) and (c)  of Proposition \ref{rel1}. 
Using the first assertion and Equations \eqref{eq1} and \eqref{eq2} we have
\[
    0 = [\sigma_{2^p, 2^{q}}] = \sum_{\ell=0}^{2^p -1}\sum_{\nu=0}^{2^q}  y_{2^p,\ell}\cdot x_{2^{q},\nu}\cdot e^{\ell+\nu}\rho^{1,1}_{2^p+\ell,2^q+ \nu} = \rho^{1,1}_{2^p, 2^{q}} + e^{2^{q}} \rho^{1,1}_{2^p, 2^{q+1}} 
\]
for all $p, q \geq 0$, where we use $y_{2^p,0} = x_{2^q, 0} = x_{2^q, 2^q} = 1$, the last equation by Lemma \ref{lem:miracle}. 
Hence we have 
\[
  \rho^{1,1}_{i,j} + e^j \rho^{1,1}_{i,2j} = 0
 \]
 if $i$ and $j$ are powers of $2$. The remaining claim follows by interchanging $i$ and $j$. 
 \end{proof} 

Finally we get the following uniform Euler torsion estimate. Here we use the relation $[\tau_m] = 0$ for $m \geq 2$. 

\begin{prop} \label{taus} We have 
\[
  e \cdot  \rho^{1,1}_{1,j} = 0  
\]
for all $j \geq 2$. 
\end{prop} 

\begin{proof} The assertion follows from Proposition \ref{sum},  if $j$ is not a power of $2$. 
It hence remains to handle the case when $j \geq 2$ is a power of $2$. 

First we need some preparation. Write the coproduct $\Delta(z) = \Delta(y(\epsilon))$ 
as in Equation \eqref{introbeta} in the form 
\[
     \Delta(y(\epsilon)) =  \sum_{i,j \geq 0 } \beta_{i,j}^{0,m} \cdot y(V_0^i) \otimes y(V_m^j)  \, , 
\]
where we henceforth assume  $m \geq 1$. For all $j  \geq 2$ we then have (recalling $\beta_{i,j}^{1,m} = \rho_{i,j}^{1,m}$ for $i+j \geq 2$) 
\begin{equation} 
  \label{eq3} \sum_{i \geq 0} \beta_{i,j}^{0,m} \cdot y(V_0^i) \otimes y(V_m^j) = 
      \sum_{i \geq 0} \beta_{i,j}^{1,m} \cdot y(V_1^i) \otimes y(V_m^j) = 
      \sum_{i \geq 0} \rho_{i,j}^{1,m} \cdot y(V_1^i) \otimes y(V_m^j) \, . 
 \end{equation} 
 According to the base change formula preceding Lemma \ref{lem:base_change} we have 
 \[
    \rho_{i,j}^{1,m} \cdot y(V_1^i) \otimes y(V_m^j)  = 
           \begin{cases} \rho_{i,j}^{1,m} \cdot y(V_0^i) \otimes y(V_m^j) \textrm{ for even } i  \\
                                                                              \rho_{i,j}^{1,m} \cdot \big( e \cdot y(V_0^{i-1}) \otimes y(V_m^j)  + y(V_0^i) \otimes y(V_m^j) \big) \textrm{ for odd } i \, , 
                                                       \end{cases} 
 \]
again using the fact that modulo the  ideal generated by the elements $\rho_{pq}$ we have $d_0 = e$, $d_1 = 1$ and $d_{i} = 0$ for $i > 1$. Comparing coefficients of the left and right hand side in Equation \eqref{eq3} we obtain 
\[
   \beta_{0,j}^{0,m} = \rho_{0,j}^{1,m} +  e \rho_{1,j}^{1,m} 
\]
for all $j  \geq 2$. We have $\rho^{1,m}_{0,j} = 0$ for $j  \geq 0$ by the co-unitality of the coproduct $\Delta$, hence 
the last equation implies 
 \[
    \beta^{0,m}_{0,j} =   e\rho^{1,m}_{1,j} 
 \]
 for all $j \geq 2$. 
 
After these preparations let $j \geq 2$ be a power of $2$. Since
$e  \rho^{1,1}_{1,j'} = 0$,  if $j' \geq 3$ is not a power of $2$ (by Proposition \ref{sum}), Equation \eqref{eq1} shows 
 \[
  \beta^{0, j}_{0,j} = e \cdot \rho^{1,j}_{1,j} = e 
  \cdot \sum_{ \ell = 0}^{j-1} y_{j,\ell} \cdot e^{\ell} \cdot \rho^{1, 1}_{1, j+\ell}   = e \cdot \rho^{1,1}_{1,j} \, , 
\] 
where the first equation follows from the preceding remarks. We therefore get 
\[
   0 = [\tau_{j}] = \beta^{0,j}_{0,j} = e \cdot \rho^{1,1}_{1,j}  
\]
as required. This finishes the proof of Proposition \ref{taus}.
\end{proof} 

Now let $i = 2^p$, $j = 2^q$, $p, q \geq 0$,  where we assume $p \leq q$ without loss of 
generality. Applying Proposition \ref{sum} several times and using Proposition \ref{taus} we get 
\[
    \rho^{1,1}_{i,j} = e^j \rho^{1,1}_{2^p,2j} =  e^{j-2^{p-1} } \rho^{1,1}_{2^{p-1} , 2j} = \cdots = e^{j - ( 2^{p} -1)} \rho^{1,1}_{1,2j} = 0 .
 \] 
This finishes the proof of Theorem \ref{step2}.

\end{document}